\documentclass{amsart}
\usepackage[all]{xy}

\usepackage{amsmath,amssymb,calligra,mathrsfs}
\usepackage[all]{xy}

\newtheorem{De}{Definition}

\newtheorem{Le}[De]{Lemma}

\def\t{\otimes }

\def \HL{{\sf HL}}
\def \H{{\sf H}}
\def \g{{\frak g}}

\def \r{{\frak r}}
\def \s{{\frak s}}

\def\t{\otimes }

\begin{document}

\title{A counterexample to a Proposition of Feldvoss-Wagemann and Burde-Wagemann}

\author[T. Pirashvili]{T. Pirashvili}
%\address{Razmadze mathematical Institute}

\maketitle
\section{Introduction}
In what follows  $\g$ is a  finite dimensional Lie algebra over the field of complex numbers. Around 1990 J.-L. Loday discovered that the the Chevalley-Eilenberg \cite{Ch}   boundary operator $d:\Lambda^{n}(\g)\to\Lambda^{n-1}(\g)$ has a canonical lift as the boundary operator $\delta:\g^{\t n}\to \g^{\t n-1}$, such that the canonical quotient map $\g^{\t n}\to \Lambda^n\g$ defines a morphism of chain complexes 
$$\xymatrix{\cdots  \ar[r]^\delta  &  g^{\t 3}\ar[d]  \ar[r]^\delta    & \g^{\t 2}    \ar[r]^\delta \ar[d]&\g \ar[d]^{id}\\\ 
\cdots   \ar[r]^d & \Lambda^3(\g)   \ar[r]^d & \Lambda^2 (\g) \ar[r]^d &\g}$$
Here   $d$  and  $\delta$ are given by
$$d(x_1\wedge \cdots \wedge x_n)=\sum_{1\leq i<j\leq n}(-1)^{i+j-1}[x_i,x_j] \wedge x_1\wedge \cdots \hat{x}_i \wedge \cdots \wedge \hat{x}_j \wedge \cdots \wedge x_n.  
$$
$$\delta(x_1\t\cdots \t x_n)=\sum_{1\leq i<j\leq n}(-1)^{j}x_1\t \cdots \t x_{i-1} \t [x_i,x_j]\t \cdots \t \hat{x}_j\t \cdots \t  x_n$$
The homology of the top complex is known as the Leibniz homology of $\g$ and they are denoted by $\HL_*(\g)$, while the homology of the bottom complex is the homology of $\g$ and is denoted by $\H_*(\g)$. In fact, Leibniz homology is defined in a much wider class of algebras known as Leibniz algebras, see \cite{HC}, \cite{UL} and \cite{LR}

In 1994 Ntolo \cite{ntolo} and the author \cite{HL} independently proved that if  $\g$ is semi-simple then $\HL_i(\g)=0$ for all $i>0$. Whether the converse is also true is our  
weak conjecture \cite{perfect}. The strong conjecture of \cite{perfect}  says that $\g$ is a semi-simple iff $\HL_i(\g)=0$ for all $i>>0$ \cite{perfect}. 

In a recent preprint \cite{bw} the authors claimed that they proved the weak-conjecture. The aim of this note is to construct  a counterexample to the crucial Proposition 3.1 of \cite{bw}.  This also  imply that  at least Theorem 4.7 and Corollary 4.8 of \cite{left_lb} are wrong.

Thus the weak conjecture of \cite{perfect} is still open.

\section{Notations}
By abuse of notation, the adjoint representation of $\g$ is also denoted by $\g$.  Denote by $\r$ the solvable radical of $\g$ and put $\s=\g/\r$.  So, $\s$ is either trivial or semi-simple Lie algebra.  In any case,  one has a split short exact sequence of Lie algebras
$$0\to \r\to\g\to \s\to 0.$$
By choosing such a section,  we obtain an action of $\s$ on $\r$. By functoriality $\s$ act also  on $\H_*(\r)$ and $\H_*(\r,\r)$.
In general these actions are highly nontrivial. 
%Here $\H_*(\g,M)$ (resp. $\H^*(\g,M)$) denotes the Chevalley-Eilenberg (co)cohomology of a Lie algebra $\g$ with coefficients in a representation $M$. If $M=\mathbb{C}$ with trivial action of $\g$, we simply write $\H_*(\g)$ (resp. $\H^*(\g)$). We willl also need Leibniz (co)homology  $\HL_*(\g,M)$ (resp. $\HL^*(\g,M)$), first introduced by Loday in \cite{HC} and further study in \cite{UL}, \cite{HL}.The aim of this work is to give a counterexample to .  As  a consequence of this article 

For a vector space $V$, we let $V^\sharp$ the dual of $V$.

\begin{Le}\label{1} For a Lie algebra $\g$ the following conditions are equivalent:

(i) The restriction map $\g^\sharp \to \r^\sharp$ yields an isomorphism
$$\HL^*(\g,\g^\sharp)\to \HL^*(\g,\r^\sharp)$$

(ii) The canonical map $\HL_*(\g,\r)\to \HL_*(\g,\g)$ is an isomorphism

(iii) One has $\HL_*(\g,\s)=0$.

(iv) One has $\H_*(\g,\s)=0$.

(v) One has $\H_0(\s, \H_*(\r)\t \s))=0.$

\end{Le}

\begin{proof}  The equivalence (i) $ \Longleftrightarrow$ (ii) follows from the fact that the cohomologies are dual vector spaces of homologies. 

The equivalence (ii)  $ \Longleftrightarrow$ (iii) follows from the homological long exact sequence associated to the short exact sequence  of $\g$-modules $0\to \r \to \g \to \s\to 0$.

We now  show  (iii)  $ \Longrightarrow$ (iv). Assume $\HL_*(\g,\s)=0$. Then by   Theorem A of \cite{HL} we have $\H^{rel}_*(\g,\s)=0$ and thus by Proposition 1 of \cite{HL} we have $\HL_i(\g,\s)=\H_i(\g,\s)$ for $i\geq 2$. Since these groups are also isomorphic for $i=0,1$ (this a general fact), we obtain $\H_*(\g,\s)=0$.

To show  (iv)  $ \Longrightarrow$ (iii), assume  $\H_*(\g,\s)=0$. We will prove by induction on $i$ that $\HL_n(\g,\s)=0$. This is obvious for $n=0,1$ because in these dimensions Leibniz and Lie homology are the same. Assume we have proved that $\HL_i(\g,\s)=0$ for all $0\leq i\leq  n.$ By Theorem $A$ \cite{HL} we obtain $\HL_i^{rel}(\g,\s)=0$ for all $0\leq i\leq n$. Hence by Proposition 1 \cite{HL} the canonical map $\HL_i(\g,\s)\to \H_i(\g,\s)$ is an isomorphism for all $0\leq i\leq n+2$ and the induction step works. Thus $\HL_*(\g,\s)=0$.

To show  (iv)  $ \Longleftrightarrow$ (v),  We consider the Lie algebra extension $0\to \r \to \g \to \s\to 0$. Since $\s$ is semi-simple or zero, we can use the homological version of Theorem 13 of \cite{hs} to obtain
$$\H_{n}(\g,\s)=\bigoplus_{p+q=n}\H_p(\s)\t \H_0(\s,\H_q(\r,\s))$$
Since the action of $\r$ on $\s$ is trivial, we can rewrite
$$\H_{n}(\g,\s)=\bigoplus_{p+q=n}\H_p(\s)\t \H_0(\s,\H_q(\r) \t \s)).$$
This show that v) imply iv). Since $\H_0(\s)$ one-dimensional, we see that the condition iv) also imply v). 
\end{proof}
\section{A counterexample}
Proposition 3.1 of \cite{bw} claims that for any $\g$ and $p\geq 1$ one has
$$\HL^p(\g)\cong \HL^{p-1}(\g, \r^\sharp).$$
We will introduce an example, which shows that  this holds not always.

Since $\HL^{*+1}(\g)\cong \HL^{*}(\g, \g^\sharp)$ holds always,  Proposition 3.1 of \cite{bw} is equivalent to the claim that the equivalent conditions of Lemma \ref{1} holds for all $\g$. Now we consider the six dimensional Lie algebra, for which this is not the case.  

We take  ${\s}= {\bf sl}_2$ and $\r$ to be the abelian Lie algebra, which as a module over $\s$ is the adjoint representation. Thus $\g$ is the semi-direct product $\r \rtimes \s$. Since  $\H_1(\r)=\s$, the Killing form defines a nontrivial $\s$-module map $\H_1(\r)\t \s\to \mathbb{C}.$ Hence $\g$ does not satisfy the condition v).

\end{document}